\documentclass[12pt]{article}
\usepackage[final]{epsfig}
\usepackage{graphics}
\usepackage{amsmath}
\usepackage{amsfonts}
\usepackage{latexsym}
\usepackage{amssymb}
\usepackage{graphicx}

\newtheorem{lemma}{Lemma}[section]
\newtheorem{proposition}[lemma]{Proposition}
\newtheorem{remark}[lemma]{Remark}
\newtheorem{example}[lemma]{Example}
\newtheorem{theorem}{Theorem}

\newtheorem{corollary}[lemma]{Corollary}

\newtheorem{conjecture}[theorem]{Conjecture}

\begin{document}
\newcommand{\eps}{{\varepsilon}}
\newcommand{\g}{{\gamma}}
\newcommand{\G}{{\Gamma}}
\newcommand{\proofend}{$\Box$\bigskip}
\newcommand{\C}{{\mathbf C}}
\newcommand{\Q}{{\mathbf Q}}
\newcommand{\R}{{\mathbf R}}
\newcommand{\Z}{{\mathbf Z}}
\newcommand{\RP}{{\mathbf {RP}}}
\newcommand{\CP}{{\mathbf {CP}}}
\newcommand{\Tr}{{\rm Tr\ }}
\def\proof{\paragraph{Proof.}}

\title {Variations on R. Schwartz's inequality for the Schwarzian derivative}
\author{Serge Tabachnikov\\
{\it Department of Mathematics, Penn State}\\
{\it University Park, PA 16802, USA}\\
e-mail: {\it tabachni@math.psu.edu}}
\date{\today}
\maketitle
\begin{abstract}
\noindent R. Schwartz's inequality provides an upper bound for the Schwarzian derivative of a parameterization of a circle in the complex plane and on the potential of Hill's equation with coexisting periodic solutions. We prove a discrete version of this inequality and obtain a version of the planar Blaschke-Santalo inequality for not necessarily convex polygons. We consider a centro-affine analog of L\"uk\H{o}'s inequality for the average squared length of a chord subtending a fixed arc length of a curve -- the role of the squared length played by the area -- and prove that the central ellipses are local minima of the respective functionals on the space of star-shaped centrally symmetric curves. We conjecture that the central ellipses are global minima. In an appendix, we relate the Blaschke-Santalo and Mahler inequalities with the asymptotic dynamics of outer billiards at infinity.
\end{abstract}

\section{Introduction and statements of results} \label{intro}

{\it Hill's equation} $f''(t)+k(t)f(t)=0$ is closely related with 1-dimensional projective and 2-dimensional centro-affine differential geometry. If $x(t)$ and $y(t)$ are two linearly independent solutions of the Hill equation then the ratio $y(t):x(t)$ gives a map $\R\to \RP^1$, a non-degenerate parametric curve in the projective line, and a different choice of solutions gives a projectively-equivalent curve. This provides a one-to-one correspondence between projective equivalence classes of non-degenerate curves in $\RP^1$ and second order differential operators $d^2/dt^2+k(t)$. See \cite{OT1} for basics of 1-dimensional projective differential geometry and Hill's equation. 

One may lift this solution curve from $\RP^1$ to a star-shaped curve $\gamma (t)$ in the plane satisfying the same equation 
\begin{equation} \label{Hill}
\gamma''(t)+k(t)\gamma(t)=0.
\end{equation}
 The lift is determined by the unit Wronskian condition 
 \begin{equation} \label{Wron}
 [\gamma (t),\gamma'(t)]=1,
 \end{equation}
  where $[\ ,\,]$ is the area form (that is, the determinant of two vectors). The curve $\gamma$ is  defined uniquely, up to  linear area-preserving transformations.  
 
Assume that all solutions of the Hill equation are $T$-periodic. Then the curve (\ref{Hill}) is  centrally symmetric, $\gamma(t+T)=-\gamma(t)$, and $2T$-periodic; in particular, the potential $k(t)$ is also $T$-periodic.  The quantity $T\int_0^T k(t)\ dt$ is called the {\it Lyapunov integral}, it plays an important role in the study of Hill's equation. For everywhere positive $k(t)$, the following inequality was proved in \cite{PB}:
\begin{equation} \label{Petty}
T\int_0^T k(t)\ dt \leq \pi^2,
\end{equation}
with equality only for constant $k(t)$. 
(Let us  mention in this regard a series of papers by Guggenheimer \cite{Gu1}--\cite{Gu6} on geometric theory of second-order differential equations, in particular, on Hill's equations with coexisting periodic solutions. )

Inequality (\ref{Petty}) is deduced in \cite{PB} from the 2-dimensional {\it Blaschke-Santalo inequality}.
Let $\gamma(t)$ be a smooth convex plane curve containing the origin in its interior. Fix an area form in the plane; then the dual plane also acquires an area form. The polar dual curve $\gamma^*(t)$ lies in the dual plane and consists of the covectors satisfying the two conditions
\begin{equation} \label{duality}
\gamma(t) \cdot \gamma^*(t)=1,\quad  \gamma'(t) \cdot \gamma^*(t)=0,
\end{equation}
where $\cdot$ is the pairing between vectors and covectors. The dual curve $\gamma^*$ is also convex and star-shaped. Let $A(\gamma)$ and $A(\gamma^*)$ be the areas bounded by $\gamma$ and $\gamma^*$, and assume that $\gamma$ is centrally symmetric with respect to the origin. In this case, the 2-dimensional Blaschke-Santalo inequality states that
\begin{equation} \label{BSsm}
A(\gamma) A(\gamma^*) \leq \pi^2,
\end{equation}
with equality only when $\gamma$ is a central ellipse (the Blaschke-Santalo inequality holds for not necessarily origin-symmetric convex curves; then one considers polar duality with respect to a special, Santalo, point, the point that minimizes $A(\gamma) A(\gamma^*)$). The product $A(\gamma) A(\gamma^*)$ is a centro-affine invariant of $\gamma$. See \cite{Lut1} concerning the Blaschke-Santalo and related affine geometric inequalities. 

The relation between inequalities (\ref{Petty}) and (\ref{BSsm}) is as follows. If $k>0$ then the curve $\gamma$ is convex. Use the area form $[\ ,\,]$ to identify the plane with its dual plane. Under this identification, (\ref{duality}) holds for $\gamma^*(t)=\gamma'(t)$. We have $A(\gamma)=T$ since $[\gamma,\gamma']=1$, and since $[\gamma',\gamma'']=k$, one has: $A(\gamma^*)=\int k(t)\ dt$. Thus (\ref{Petty}) follows from (\ref{BSsm}). 

Independently of \cite{PB}, R. Schwartz \cite{Sch1}, in his study of a projectively natural flow on the space of diffeomorphisms of a circle, considered a diffeomorphism $f:\R/2\pi\Z\to S^1\subset \C$ and proved  the following {\it Average Lemma}:
\begin{equation} \label{average}
\int_0^{2\pi} S(f)\ dt \leq \pi,
\end{equation}
where 
$$
S(f)=\frac{f'''}{f'}-\frac{3}{2}\left(\frac{f''}{f'}\right)^2
$$
is the Schwarzian derivative of $f(t)$ (the Schwarzian  is real for $|f(t)|=1$). See also \cite{Sch2} where  a similar inequality for a convex curve in $\RP^2$ is proved; we do not dwell on this other inequality of R. Schwartz here.

A stereographic projection from a point of a circle identifies the circle in the complex plane with the real projective line, and $f$ can be considered as a $2\pi$-periodic curve in $\RP^1$ (a different choice of the center of stereographic projection gives a projectively equivalent curve). This curve corresponds to Hill's equation whose potential, $k(t)$, can be reconstructed as the Schwarzian derivative of the ratio of its solutions, see \cite{Fl,OT1,OT2}. A computation  reveals  that inequality (\ref{average}) has the same form as (\ref{Petty}), but without the positivity assumption $k(t)>0$.\footnote{The Schwarzian  derivative is intimately related with curvature, in spherical \cite{OT0}, Lorentz \cite{Ghy,D-O,Ta}, and  hyperbolic \cite{Si} geometries.} 

\begin{figure}[hbtp]
\centering
\includegraphics[width=4in]{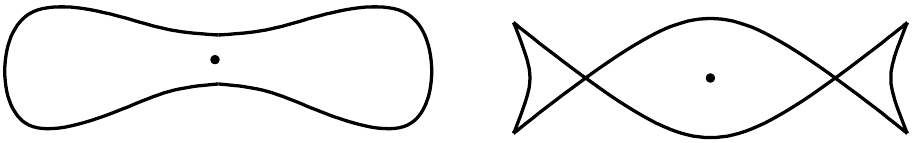}
\caption{A non-convex star-shaped curve and its polar dual}
\label{fronts}
\end{figure}

Figure \ref{fronts}, left, depicts a non-convex star-shaped curve $\gamma$. The polar dual curve $\gamma^*$ is still star-shaped in that no tangent line passes through the origin; however it has cusps, corresponding to inflections of $\gamma$ (the points at which $k=0$), and self-intersections, corresponding to double tangents of $\gamma$. Such singular curves are called {\it wave fronts}. The area $A(\gamma^*)$ is defined as the integral of the 1-form $xdy$ over the wave front $\gamma^*$,  oriented so that the tangent line turns in the positive sense. The Average Lemma of Schwartz can be interpreted as a 2-dimensional Blaschke-Santalo inequality for such star-shaped curves.

See  \cite{CHLYZ, NZ} for a version of Blaschke-Santalo inequality for not necessarily convex plane curves in terms of the support function,  \cite{LYZ} for  a version of the Blaschke-Santalo inequality for compact sets, and \cite{Leh} for a functional  Blaschke-Santalo inequality.

We provide a discretization of Schwartz's inequality. Namely, we prove a version of inequalities (\ref{Petty}), (\ref{BSsm})  and (\ref{average})  for star-shaped, but not necessarily convex, polygons. 
Consider an origin-symmetric star-shaped $2n$-gon  in the plane with vertices $V_i$ in their cyclic order about the origin, such that $V_{i+n}=-V_i$ and $[V_i,V_{i+1}]=1$ for all $i$. Let $c_i=[V_{i-1},V_{i+1}]$; the sequence $c_i$ is $n$-periodic. Set $F_n=\sum_{i=1}^n c_i$. Obviously, each $c_i$, and hence $F_n$, is invariant under the action of $SL(2,\R)$ on polygons.

\begin{theorem} \label{discr1}
One has:
$$
F_n \geq 2n\cos\frac{\pi}{n},
$$
with equality only for the $SL(2,\R)$-equivalence class of regular polygons.
\end{theorem}

\begin{figure}[hbtp]
\centering
\includegraphics[width=4in]{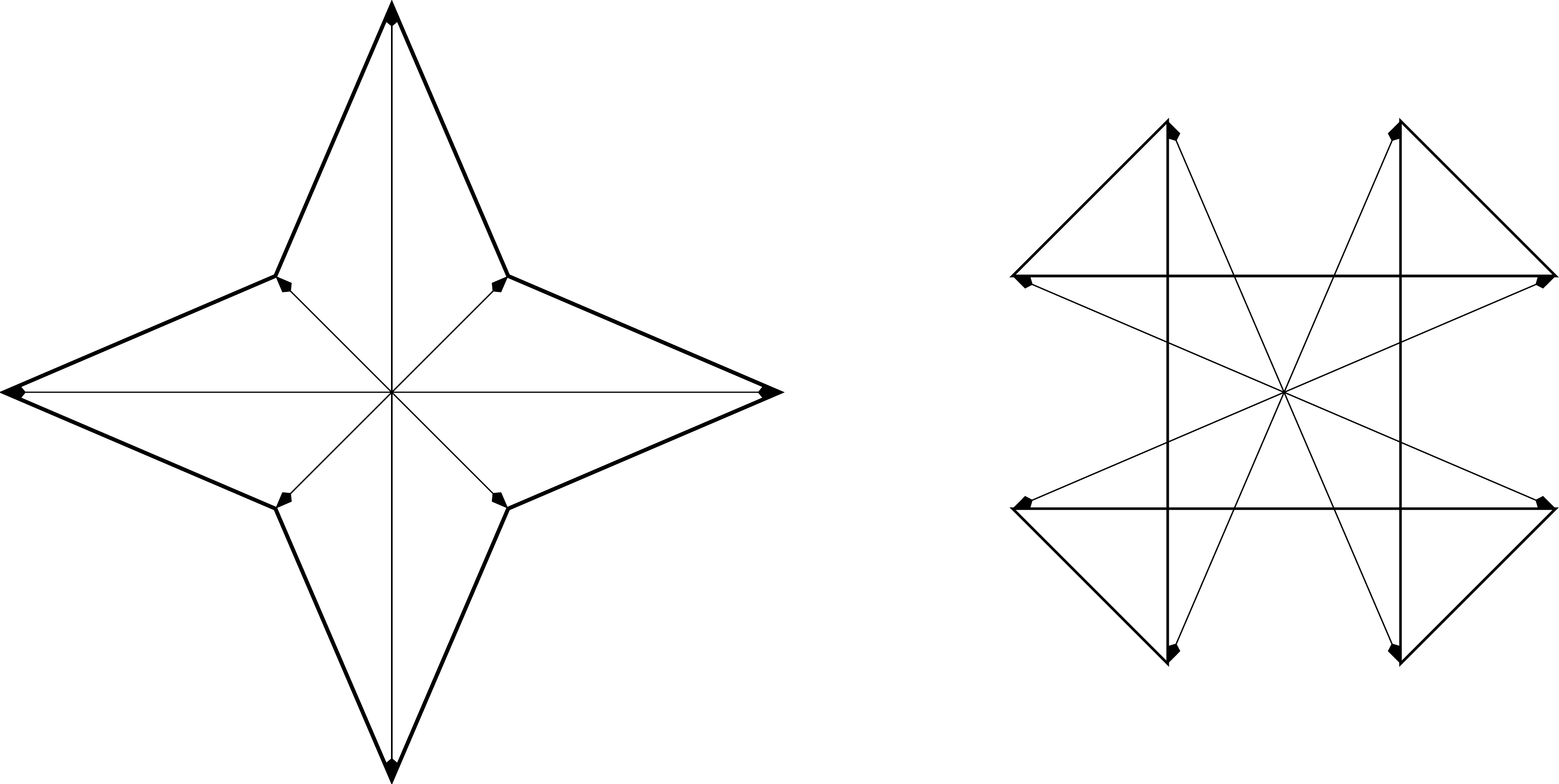}
\caption{A star-shaped polygon $V$ and its dual  $V^*$}
\label{polygons}
\end{figure}

Consider an origin-symmetric $2n$-gon $V^*$ with vertices $V^*_i =V_{i+1}-V_i$. This polygon may be self-intersecting, see Figure \ref{polygons}. We shall see that $V^*$ is polar dual to $V$. Using the same notation for areas as in (\ref{BSsm}), the following polygonal Blaschke-Santalo inequality holds.

\begin{theorem} \label{BS}
One has: 
\begin{equation} \label{BSineq}
A(V) A(V^*) \leq 4n^2\sin^2 \frac{\pi}{2n},
\end{equation}
with equality only for centro-affine regular polygons.
\end{theorem}

Note that the limit $n\to\infty$ of the  right hand side of (\ref{BSineq}) is $\pi^2$, the right hand side of (\ref{BSsm}).

Theorem \ref{discr1} is reminiscent of another extremal property of regular polygons, in terms of their diagonal lengths, discovered by  G. L\"uk\H{o}
 \cite{Luk}. Let $V_i$ be an $n$-gon and $1<k<n-1$ be  fixed. Assume that $|V_i V_{i+1}|\leq C$ for some constant $C$ and $f$ is an increasing concave function. Then
$$
\frac{1}{n}\sum_{i=1}^n f(|V_i V_{i+k}|^2) \leq f\left(C^2 \sin^2 \frac{k\pi}{n}/ \sin^2 \frac{\pi}{n}\right),
$$
with equality only for regular $n$-gons. In particular, one has an upper bound on the average length of $k$-diagonals of a polygon:
$$
\frac{1}{n}\sum_{i=1}^n |V_i V_{i+k}| \leq C \sin \frac{k\pi}{n}/ \sin \frac{\pi}{n}.
$$
In the limit $n\to\infty$, one has a similar upper bound on the average chord length for smooth curves, see  \cite{ACFGH, EHL}:
\begin{equation} \label{chord}
 \frac{1}{2\pi}\int_0^{2\pi} f(|\gamma(t+c)-\gamma(t)|^2)\ dt \leq f\left(4\sin^2 \frac{c}{2}\right),
\end{equation}
with equality only for the round unit circle; here $t$ is arc length parameter and the total length of $\gamma$ is normalized to  $2\pi$. This inequality was used in \cite{ACFGH} to prove that many knot energies are uniquely minimized by round circles. 

In the spirit of Theorem \ref{discr1}, we propose to consider a centro-affine version of inequality (\ref{chord}). Let $\gamma(t)$ be an origin-symmetric star-shaped $2\pi$-periodic curve such that $\gamma(t+\pi)=-\gamma(t)$, satisfying the unit Wronskian condition (\ref{Wron}), the centro-affine analog of arc length parameter. For $\alpha \in (0,\pi)$, set
$$
I(\alpha)=\frac{1}{\pi} \int_0^{\pi} [\gamma (t),\gamma (t+\alpha)]\ dt.
$$

\begin{conjecture} \label{conj}
For every $\alpha$, one has: 
\begin{equation} \label{conjineq}
I(\alpha) \geq \sin \alpha,
\end{equation}
with equality only for central ellipses. 
\end{conjecture}

For infinitesimal $\alpha$, the Taylor expansion up to third order shows that (\ref{conjineq}) implies inequality (\ref{Petty}) (with $T=\pi$). Thus Conjecture \ref{conj} is indeed a generalization of Schwartz's Average Lemma (\ref{average}).

In a sense, one can solve the centro-affine parameterization equation (\ref{Wron}). Let $t$ be the angular coordinate in $\R^2$; then $t$ gives a parameterization of $\RP^1$ such that $t$ and $t+\pi$ correspond to the same point. Let $f(t)$ be an orientation preserving diffeomorphism of $\RP^1$ which we consider as a diffeomorphism $f:\R\to \R$ satisfying $f(t+\pi)=f(t)+\pi$. Then the curve
\begin{equation} \label{param}
\gamma(t)=\frac{1}{\sqrt{f'(t)}} (\cos f(t),\sin f(t))
\end{equation}
satisfies (\ref{Wron}), and all solutions are obtained this way; see, e.g., \cite{OT0}. Conjecture \ref{conj} can be reformulated as follows:
\begin{equation} \label{intineq}
\frac{1}{\pi} \int_0^{\pi} \frac{\sin (f(t+\alpha)-f(t))}{\sqrt{f'(t+\alpha)f'(t)}}\ dt \geq \sin \alpha
\end{equation}
for all diffeomorphisms $f$ as above and every $\alpha \in (0,\pi)$, with equality only for projective diffeomorphisms of $\RP^1$.

We prove a weak version of Conjecture \ref{conj}.

\begin{theorem} \label{local}
For every $\alpha \in (0,\pi)$, the central ellipses are local minima of the functional $I(\alpha)$.
\end{theorem}

Namely we shall show that the central ellipses form a critical 3-dimensional manifold of $I(\alpha)$ with a  Hessian, positive definite in the normal direction.

In the spirit of  \cite{ACFGH, EHL}, one may consider the {\it areal energy} of a centro-affine parameterized curve $\gamma(t)$, as considered above:
$$
G[\gamma]=\int_0^{\pi} \int_0^{\pi} g([\gamma(t),\gamma(t+\alpha)],\alpha)\ dt\ d\alpha,
$$
where $g$ is a function of two variables. One conjectures that, for a broad class of functions $g$, this areal energy $G[\gamma]$ is uniquely minimized by the central ellipses. 

The content of the paper is as follows. In Section \ref{proofs12} we prove Theorems \ref{discr1} and \ref{BS}. Our proof of Theorem \ref{discr1} is by way of Morse theory on the space of equivalence classes of relevant polygons. Describing this space, we use some combinatorial formulas known in the theory of frieze patterns. In Section \ref{proof4}, using Fourier expansions of periodic functions, we prove Theorem \ref{local}. The proof reduces to an infinite series of trigonometric inequalities. 

Section \ref{outerbill} is an appendix devoted to a somewhat unexpected appearance of the Blaschke-Santalo and Mahler inequalities, as well as the isoperimetric inequality in plane Minkowski geometry, in the study of outer billiards, a geometrically natural dynamical system akin to the more familiar, inner, billiards. To avoid expanding this introduction any further, we postpone the discussion of outer billiards until Section \ref{outerbill}. 

\section{Proofs of Theorems \ref{discr1} and \ref{BS}} \label{proofs12}

Denote by ${\cal P}_n$ the space of origin-symmetric star-shaped $2n$-gons $(V_i)$ satisfying $V_{i+n}=-V_i$ and $[V_i,V_{i+1}]=1$ for all $i$; and let ${\cal M}_n$ be its quotient space  by $SL(2,\R)$. 

\begin{lemma} \label{mfld}
The spaces ${\cal P}_n$ and ${\cal M}_n$ are  smooth $n$- and $(n-3)$-dimensional manifolds, respectively.
\end{lemma}

\proof Consider ${\cal P}_n$ as a subvariety in $(\R^2)^n$ defined by the conditions $\varphi_i=1,\ i=0,\dots,n-1,$ where $\varphi_i=
[V_i,V_{i+1}]$.
Let $V=(V_i)$ be a polygon in ${\cal P}_n$. We want to show that our condition define a smooth submanifold: if $\sum_i \lambda_i d\varphi =0$ at $V$ then all $\lambda_i=0$. Consider  a test tangent vector $\xi=(0,\dots,0,\xi_k,0,\dots,0)\in (\R^2)^n$ where $\xi_k\in\R^2$ is at $k$-th position. Then
$$
0=\sum_i \lambda_i d\varphi (\xi)=[\xi_k, \lambda_k V_{k+1}-\lambda_{k-1} V_{k-1}].
$$
Hence $\lambda_k V_{k+1}=\lambda_{k-1} V_{k-1}$. However, the vectors $V_{k-1}$ and $V_{k+1}$ are linearly independent, so $\lambda_k=\lambda_{k-1}=0$. This holds for all $k$ establishing the claim.

Since $SL(2,\R)$ acts freely on ${\cal P}_n$, the quotient space ${\cal M}_n$ is an $(n-3)$-dimensional manifold.
\proofend

\begin{remark} \label{conf}
Polygons in the projective line and in the affine plane.
{\rm
One has a natural map from ${\cal P}_n$ to ${\cal C}_n$, the configuration space of $n$ points in $\RP^1$. If $n$ is odd, this projection is a bijection on the connected component consisting of $n$-gons with winding number 1. 

Indeed, let $v_0,\dots,v_{n-1}\in\RP^1$ be such that the segments $[v_i,v_{i+1})$ (not containing other points $v_j$) cover the projective line once. Lift  points $v_i$ to vectors $U_i\in\R^2$ so that $[U_i,U_j]>0$ for $0\leq i<j\leq n-1$. We want to rescale these vectors, $V_i=t_iU_i$, so that $[V_i,V_{i+1}]=1$ for $i=0,\dots,n-2$, and $[V_{n-1},-V_0]=1$. This gives the system of equations
$$
t_it_{i+1}=1/[U_i,U_{i+1}],\ \ i=0,\dots,n-2;\ \ t_{n-1}t_0=1/[U_0,U_{n-1}],
$$
that has a unique solution for odd $n$. This provides an inverse map ${\cal C}_n\to {\cal P}_n$.

However, if $n$ is even, the projection ${\cal P}_n\to {\cal C}_n$ has a 1-dimensional fiber given by the scaling:
$$
V_{2i}\mapsto tV_{2i},\ V_{2i+1}\mapsto t^{-1} V_{2i+1},\quad t\in\R_+.
$$ 
The image of the projection ${\cal P}_n\to {\cal C}_n$ has codimension 1; it is given by the condition 
$$
\prod_{i\ {\rm even}} [U_i,U_{i+1}] = \prod_{i\ {\rm odd}} [U_i,U_{i+1}],
$$
that does not depend on the lifting.
}
\end{remark}

We interpret the cross-products $c_i$  as follows. One can express each next vector $V_{i+1}$ as a linear combination  the previous two, and the conditions $[V_{i-1},V_i]=[V_i,V_{i+1}]=1, [V_{i-1},V_{i+1}]=c_i$ imply that the coefficients are as follows:
\begin{equation} \label{rec}
V_{i+1}=c_i V_i - V_{i-1}.\footnote{This difference equation is a discrete analog of Hill's equation.}
\end{equation}
Set: $F_{i,j}=[V_i,V_j]$.

\begin{lemma} \label{MatrixLemma}
For $j-i\geq 2$, one has:
\begin{equation} \label{matrix}
F_{i,j}=\left|\begin{array}{cccccc}
c_{i+1}&1&0&0&\dots&0\\
1&c_{i+2}&1&0&\dots&0\\
0&1&c_{i+3}&1&\dots&0\\
\dots&\dots&\dots&\dots&\dots&\dots\\
\dots&\dots&\dots&0&1&c_{j-1}
\end{array}\right|.
\end{equation}
One also has:
\begin{equation} \label{Vform}
V_j=-F_{0,j} V_{-1}+F_{-1,j} V_0.
\end{equation}
\end{lemma}

\proof Equation (\ref{matrix}) is proved by induction on $j-i$. The determinants satisfy the recurrence
$$
F_{i,j+1}=c_jF_{i,j}-F_{i,j-1},
$$
but, due to (\ref{rec}),  the same recurrence holds for the cross-products:
$$
[V_i,V_{j+1}]=c_j[V_i,V_j]-[V_i,V_{j-1}],
$$
which makes it possible to use induction.

Equation (\ref{Vform}) follows from the fact that both sides have the same cross-products with $V_{-1}$ and $V_0$.
\proofend

\begin{corollary} \label{closed}
One has:
\begin{equation} \label{mono}
F_{0,n-1}=1,\ F_{-1,n-1}=0,\ F_{0,n}=0.
\end{equation}
\end{corollary}

\proof
Since $V_{n-1}=-V_{-1}$ and $V_n=-V_0$, equation (\ref{mono}) follows from (\ref{Vform}) for $j=n-1$ and $j=n$. (There is a fourth condition, $F_{-1,n}=-1$, but it follows from the fact that the monodromy map $(V_{-1},V_0) \mapsto (V_{n-1},V_n)$ is area-preserving).
\proofend

\begin{remark}  \label{frieze}
Frieze patterns.
{\rm The cross-products $F_{i,j}$ with fixed $j-i$ form the rows of a frieze pattern; this means that the following equality holds: 
\begin{equation} \label{Fs}
F_{i-1,j-1}F_{i,j} - F_{i,j-1}F_{i-1,j}=1
\end{equation}
(see \cite{Cox} concerning frieze patterns). Indeed, if
$$
V_{j-1}=aV_{i-1}+bV_i,\ V_j=cV_{i-1}+dV_i,
$$
then, using $[V_{j-1},V_j]=1$, one has: $ad-bc=1$. On the other hand,  the left hand side of (\ref{Fs}) is equal to $-bc+ad$, and (\ref{Fs}) follows. 

Thus the numbers $c_i$ form the first non-trivial row of a frieze pattern (after a row of 0s and a row of 1s). Formula (\ref{matrix}) can be found in \cite{Cox} too. 
}
\end{remark}

The functions $c_i$ serve as coordinates in ${\cal M}_n$. They are not independent: they satisfy the three relations of Corollary \ref{closed}. One can use formula (\ref{Vform}) to reconstruct the equivalence class of a polygon from $c_i$; this is used in the next lemma. But first consider the examples of $n=4$ and $n=5$.

\begin{example} \label{n45}
{\rm $n=4$: applying a transformation from $SL(2,\R)$, we may assume that  $V_0=(1,0), V_3=(0,1)$. Then $V_1=(x,1), V_2=(1,y)$. The condition $[V_1,V_2]=1$ yields $xy-1=1$. Thus ${\cal M}_4$ is the hyperbola $xy=2, x>0$. The cross-products $c_i$ are as follows:
$c_0=x, c_1=y, c_2=x, c_3=y,$
and $F_4=2(x+y)$, which has the minimum for $x=y=\sqrt{2}=2\cos (\pi/4)$. 

$n=5$: once again, assume that $V_0=(1,0), V_4=(0,1)$. Let $V_1=(x,1), V_3=(1,y)$, and $V_2=(a,b)$. Then the conditions $[V_1,V_2]=[V_2,V_3]=1$ yield $x=(1+a)/b, y=(1+b)/a$. The cross-products $c_i$ are as follows:
$$
c_0=x=\frac{1+a}{b}, c_1=b, c_2=xy-1=\frac{1+a+b}{ab}, c_3=a, c_4=y=\frac{1+b}{a},
$$
and 
$$
F_5=a+b+\frac{1+a}{b}+\frac{1+b}{a}+\frac{1+a+b}{ab}.
$$
The only critical point of this function is $a=b=(1+\sqrt{5})/2$, the golden ratio; then $c_i=2\cos(\pi/5)$ for all $i$.
}
\end{example} 

Consider the function $F_n: {\cal M}_n \to \R$.

\begin{lemma} \label{proper}
$F_n$ is a proper function, that is, the set $S_m:=\{F_n\leq m\}$ is compact for every positive constant $m$.
\end{lemma}

\proof Consider a sequence of polygons $V^j=(V_i^j),\ j=1,2,\dots$  in $S_m$ and let $c_i^j$ be the respective sequence of the cross-products $[V_{i-1}^j,V_{i+1}^j]$. Since $c_i^j>0$, we have $c_i^j\leq m$, and hence, considering a subsequence if necessary, we may assume that $c_i^j \to \bar c_i$ as $j \to \infty$. Since $F_{i-1,i+2}^j >0$, it follows from (\ref{matrix}) that $\bar c_i \bar c_{i+1} \geq 1$. Therefore, all $\bar c_i$ are separated from 0.

We can use the numbers $\bar c_i$ to construct a polygon in ${\cal M}_n$ which is the limit of the (sub)sequence $V^j$. Namely, choose two vectors, $\bar V_{-1}$ and $\bar V_0$ with $[\bar V_{-1},\bar V_0]=1$, and use the ``barred" version of recurrence  (\ref{Vform}) to construct a polygon. 
The periodicity condition $\bar V_{n-1}=-\bar V_{-1}, \bar V_n=-\bar V_0$ follows from the fact that equations (\ref{mono}) still hold in the limit. The resulting polygon is star-shaped because, in the limit, $F_{i,j} \geq 0$ for $0\leq i<j\leq n-1$ as well.
\proofend

Now we describe the critical points of the function $F_n$ in ${\cal M}_n$.

\begin{lemma} \label{crit}
A polygon $V=(V_i)$ is a critical point of the function $F_n$ if and only if all the cross-ratios are equal,
$c_i=2\cos(\pi/n)$, and the polygon is affine-regular.
\end{lemma}

\proof
Assume that $n\geq 6$ (otherwise, see Example \ref{n45}). Consider six consecutive vertices $V_{i-2},\dots,V_{i+3}$. Consider an infinitesimal deformation 
\begin{equation} \label{deform}
V_i \mapsto V_i + \varepsilon V_{i-1},\ V_{i+1} \mapsto V_{i+1} + \delta V_{i+2}.
\end{equation} 
This deformation does not change the cross-products $[V_{i-1},V_i]$ and $[V_{i+1},V_{i+2}]$. For $[V_i,V_{i+1}]$ to remain the same in the linear approximation, one needs to have $\varepsilon[V_{i-1},V_{i+1}] + \delta [V_i,V_{i+2}]=0$. Hence $\varepsilon=t c_{i+1}, \delta =-tc_i$ where $t$ is an infinitesimal. In particular, $\varepsilon+\delta=t(c_{i+1}-c_i)$.

Next one computes the rate of change of $F_n$ under the deformation (\ref{deform}). This equals
\begin{equation} \label{change}
\varepsilon [V_{i-2},V_{i-1}]+\delta [V_{i-1},V_{i+2}] + \varepsilon [V_{i-1},V_{i+2}] + \delta [V_{i+2},V_{i+3}] = (\varepsilon+\delta)(F_{i-1,i+2}+1).
\end{equation}
By (\ref{matrix}), $F_{i-1,i+2}=c_ic_{i+1}-1$, hence (\ref{change}) equals
$(\varepsilon+\delta)c_ic_{i+1} = t(c_{i+1}-c_i)c_ic_{i+1}.$
This is zero if and only if $c_i=c_{i+1}$. Therefore, if a point is critical then all $c_i$ are equal. This is the case of the (affine) regular polygon.

It remains to check that (an equivalence class of) the regular polygon $V=(V_i)$ is a critical point of $F_n$. Consider an infinitesimal deformation $V_i\mapsto V_i + \varepsilon U_i$ where $U_i=a_iV_{i-1}+b_iV_{i+1}$. Since the deformation does not change $[V_i,V_{i+1}]$, one has: $[U_i,V_{i+1}]+[V_i,U_{i+1}]=0$, that is, $a_i c_i + b_{i+1}c_{i+1}=0$. Since all $c_i$ are equal, $c_i=c$ for all $i$, one has $a_i=-b_{i+1}$, and in particular, $\sum(a_i+b_i)=0$. Finally,  the rate of change of $F_n$ equals 
\begin{equation*}
\begin{split}
\sum [U_{i-1},V_{i+1}]&+[V_{i-1},U_{i+1}]=\\
a_{i-1} [V_{i-2},V_{i+1}]+b_{i-1}[V_{i},V_{i+1}]&+a_{i+1}[V_{i-1},V_{i}]+b_{i+1}[V_{i-1},V_{i+2}] =\\ 
a_{i-1} (c_{i-1}c_i-1)+b_{i-1}+a_{i+1}&+b_{i+1}(c_ic_{i+1}-1) = c^2 \sum(a_i+b_i)  =0.
\end{split}
\end{equation*}
Thus $dF_n (V)=0$.
\proofend

Now we can prove Theorem \ref{discr1}. Fix a generic and sufficiently large constant $m$ and consider the manifold with boundary $S_m=\{F_n\leq m\}$. By Lemma \ref{proper}, $S_m$ is compact. Hence $F_n$ assumes minimum on it, and by Lemma \ref{crit}, this minimum corresponds to the affine-regular polygon. The respective value of the function $F_n$ is $2n\cos(\pi/n)$. This finishes the proof of Theorem \ref{discr1}.

 Let $V=(V_i)$ be a star-shaped polygon in the plane. The dual polygon $V^*=(V^*_i)$ in the dual plane is characterized by the equalities:
$$
V_i\cdot V^*_i=1,\quad \ker V^*_i = (V_{i+1} V_i),
$$
for all $i$.  As in Section \ref{intro}, we use the area form to identify the plane with its dual. 

\begin{lemma} \label{dual}
The dual polygon $V^*$ is given by $V^*_i=V_{i+1}-V_i$, and its signed area satisfies $A(V^*)=2n-F_n$.
\end{lemma}

\proof One has:
 $$
[V_i,V_{i+1}-V_i]=1\quad {\rm and} \quad [V_{i+1}-V_i,V_{i+1}-V_i]=0,
$$ 
as needed. Next,
$$
A(V^*)=\sum_1^{n} [V_{i}-V_{i-1}, V_{i+1}-V_i]=\sum_1^{n} (2-c_i)=2n-F_n,
$$
as claimed.
\proofend

Now we prove Theorem \ref{BS}. Consider a centrally  symmetric star-shaped $2n$-gon $V$ satisfying the assumptions of the theorem.    
One has $A(V)=n$, hence 
$$
A(V) A(V^*)=2n^2-nF_n \leq 2n^2-2n^2 \cos\frac{\pi}{n}=4n^2 \sin^2 \frac{\pi}{2n}, 
$$
where the inequality in the middle follows from Theorem \ref{discr1}. 

\section{Proof of Theorem \ref{local}} \label{proof4}

First of all, we show that inequality (\ref{intineq}) is indeed equivalent to (\ref{conjineq}).

\begin{lemma} \label{eq}
For a curve $\gamma$ as in (\ref{param}), one has:
$$
I(\alpha)=\frac{1}{\pi} \int_0^{\pi} \frac{\sin (f(t+\alpha)-f(t))}{\sqrt{f'(t+\alpha)f'(t)}}\ dt.
$$
\end{lemma}

\proof
Using complex notation, one has:
$$
\gamma(t)=(f'(t))^{-1/2} e^{if(t)},\ \gamma(t+\alpha)=(f'(t+\alpha))^{-1/2} e^{if(t+\alpha)},
$$
hence
$$
[\gamma(t),\gamma(t+\alpha)]=(f'(t)f'(t+\alpha))^{-1/2} [e^{if(t)}, e^{if(t+\alpha)}],
$$
and it remains to use the fact that $[\exp(i\phi),\exp(i\psi)]=\sin (\psi-\phi)$.
\proofend

We use the formula of Lemma \ref{eq} to investigate the functional $I(\alpha)$.

\begin{lemma} \label{critical}
For each $\alpha$, the central ellipses are critical points of the functional $I(\alpha)$. 
\end{lemma}

\proof
Without loss of generality, assume that we are given a unit circle parameterized by the angle parameter $t$, that is, $f_0(t)=t$. Consider an infinitesimal perturbation $f(t)=t+\varepsilon g(t)$ where $g$ is a $\pi$-periodic function. Then
$$
f(t+\alpha)=t+\alpha+ \varepsilon g(t+\alpha),\ f'(t)=1+ \varepsilon g'(t),\ f'(t+\alpha)=1+\varepsilon g'(t+\alpha).
$$
Denoting $f(t+\alpha)$ and $g(t+\alpha)$ by $f_+$ and $g_+$ respectively, one has 
$$
\sin (f_+-f)=\sin \alpha + \varepsilon  (g_+-g) \cos \alpha,\ f' f_+'=1+ \varepsilon (g_+'+g'),
$$
and hence
$$
\frac{1}{\pi} \int_0^{\pi} \frac{\sin (f(t+\alpha)-f(t))}{\sqrt{f'(t+\alpha)f'(t)}}\ dt = \sin \alpha + \frac{\varepsilon}{\pi}  \int_0^{\pi} \left((g_+-g) \cos \alpha -\frac{1}{2} (g_+'+g') \sin \alpha\right)\ dt.
$$
The last integral vanishes because 
\begin{equation} \label{zeros}
\int_0^{\pi} g_+(t)\ dt = \int_0^{\pi} g(t)\ dt\quad {\rm and}\quad \int_0^{\pi} g_+'(t)\ dt = \int_0^{\pi} g'(t)\ dt =0,
\end{equation}
as needed.
\proofend

Next we compute the Hessian of the functional $I(\alpha)$ at function $f_0(t)=t$. Write: $f(t)=t+ \varepsilon g(t) + \varepsilon^2 h(t)$, where $g$ and $h$ are $\pi$-periodic, and use the same notation as in the proof of the preceding lemma. 

\begin{lemma} \label{Hess}
One has: 
$$
I(\alpha)=\sin \alpha + \frac{\varepsilon^2}{4} \left( \sin \alpha \int_0^{\pi} (g'g_+'+3(g')^2+4gg_+-4g^2)\ dt -4\cos \alpha  \int_0^{\pi} g_+g'\ dt\right) 
$$
where terms of order 3 and higher in $\varepsilon$ are suppressed.
\end{lemma}

\proof
The computation is similar to the previous proof, but this time, one considers Taylor expansions in $\varepsilon$ up to second order. From Lemma \ref{critical} we know that the linear term in $\varepsilon$ in the expansion of $I(\alpha)$ vanishes. The quadratic term in the integrand is as follows:
$$
\sin\alpha\left(-\frac{1}{2}(g_+-g)^2+\frac{3}{8}(g_+'+g')^2-\frac{1}{2} (h_+'+h'+g_+'g')\right)
$$
$$
+\cos\alpha \left(-\frac{1}{2}(g_+-g)(g_+'+g') +(h_+-h)\right).
$$
As before, we simplify the integrals using equations (\ref{zeros}) for function $h$, similar equations for functions $g^2$ and $g_+^2$, and integration by parts
$$
  \int_0^{\pi} g(t)g_+'(t)\ dt = - \int_0^{\pi} g'(t) g_+(t)\ dt
$$
to obtain the stated result.
\proofend

Next we consider the Fourier expansion of  the $\pi$-periodic function $g(t)$ and express the Hessian in terms of the Fourier coefficients. Let 
$$
g(t)=\sum_{n\in\Z} z_n e^{i n t},\ \ z_n\in\C,\ \ z_{-n}=\bar z_n,\ \ n\ {\rm even}.
$$

\begin{lemma} \label{Fourier}
Up to a multiplicative positive constant, the quadratic part of $I(\alpha)$, as given in Lemma \ref{Hess}, is as follows:
$$
\sum_{n>0,\ {\rm even}} |z_n|^2 [(3n^2-4)\sin\alpha+(n^2+4)\sin\alpha\cos n\alpha -4n\cos\alpha\sin n\alpha].
$$
\end{lemma}

\proof
First of all, we notice that the quadratic part of $I(\alpha)$ vanishes if $g$ is a constant. Hence we may assume that $z_0=0$. Next, we have:
$$
g'=\sum in z_n e^{i n t},\ g_+=\sum z_n e^{in\alpha} e^{i n t},\ g_+'=\sum in z_n e^{in\alpha} e^{i n t}.
$$
Now we use Lemma \ref{Hess} and the fact that 
$$
\int_0^{2\pi} e^{i n t} e^{i m t}\ dt =0,
$$
unless $m=-n$, in which case the integral equals $2\pi$ (we deal with even harmonics, hence we may take the limits in the integrals to be $0$ and $\pi$). Using this fact, each integral from Lemma \ref{Hess} can be expressed in terms of the coefficients $z_n$. Let us illustrate this for the term $g_+ g$; other cases are similar:
$$
\int_0^{\pi} g_+(t) g(t)\ dt = \frac{1}{\pi} \sum_{n\ {\rm even}} z_{-n} z_n e^{in\alpha} = \frac{1}{\pi} \sum_{n\ {\rm even}}   |z_n|^2 e^{in\alpha} =
$$
$$
 \frac{2}{\pi} \sum_{n>0,\ {\rm even}}  |z_n|^2\ \frac{e^{in\alpha}+e^{-in\alpha}}{2} = \frac{2}{\pi} \sum_{n>0,\ {\rm even}}  |z_n|^2 \cos n\alpha.
$$
Collecting terms and canceling a common positive factor yields the result.
\proofend

It remains to consider the function 
\begin{equation} 
f_n (\alpha):=(3n^2-4)\sin\alpha+(n^2+4)\sin\alpha\cos n\alpha -4n\cos\alpha\sin n\alpha.
\end{equation}
We observe that $f_0(\alpha)=0$ and $f_2(\alpha)=0$. The space of even harmonics of order $\leq 2$ is 3-dimensional. 
This corresponds to the fact that the space of central ellipses is 3-dimensional: they constitute the $SL(2,\R)$-orbit of the unit circle. 
Thus Theorem \ref{local} will be proved once we establish the following fact.

\begin{figure}[hbtp]
\centering
\includegraphics[width=2.5in]{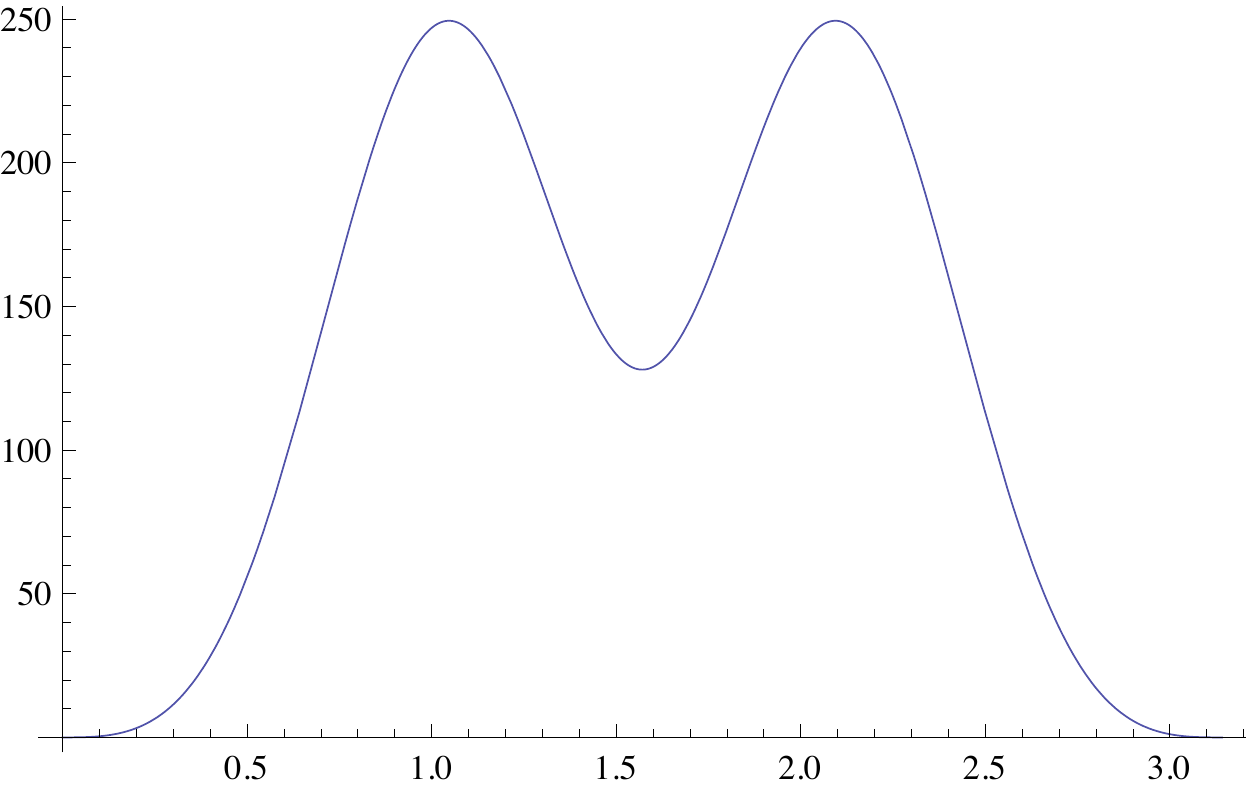}
\quad
\includegraphics[width=2.5in]{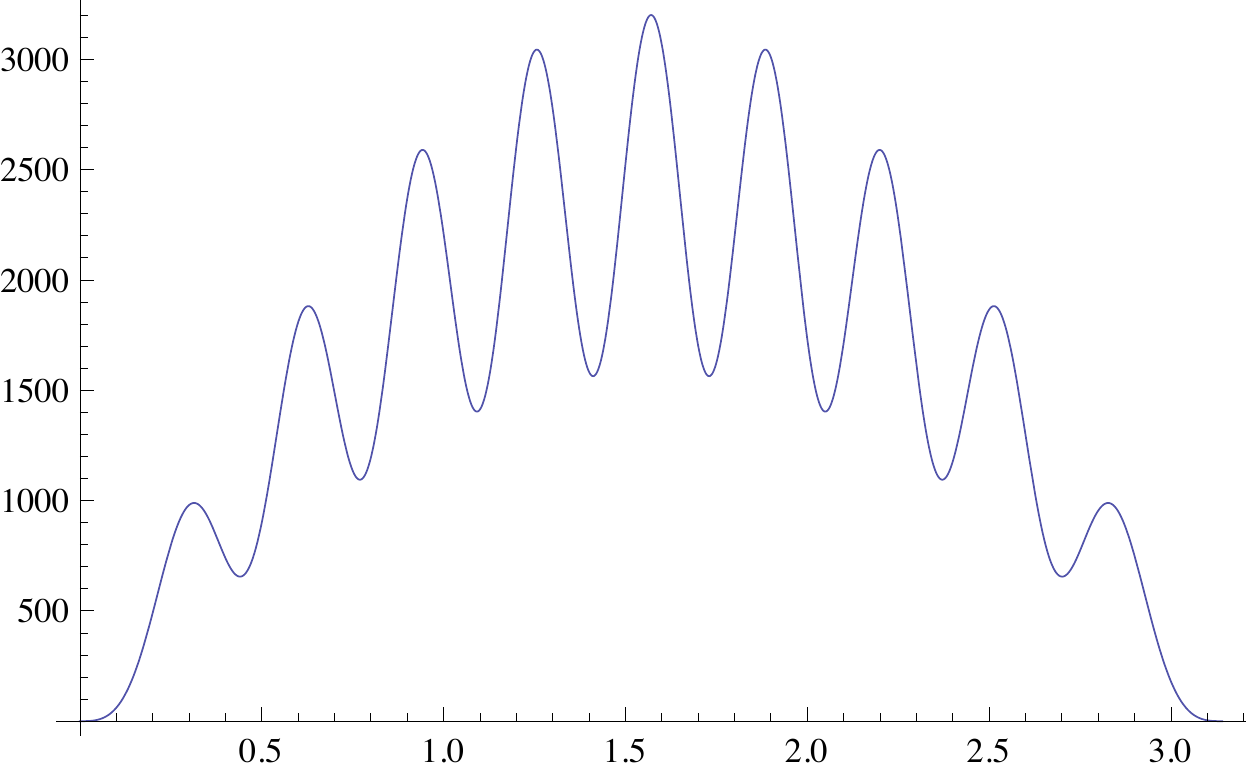}
\caption{The graphs of $f_n$ for $n=6$ and $n=20$}
\label{graphs}
\end{figure}

\begin{proposition}
For every  even $n\geq 4$ and $\alpha\in(0,\pi)$, one has: $f_n (\alpha) >0$ (see Figure \ref{graphs}).
\end{proposition} 

\proof
First, we show that $f_n(\alpha)>0$ if $\alpha$ is sufficiently separated from $0$ and $\pi$. In fact, by symmetry, we may assume that $\alpha\in(0,\pi/2)$. 
One has an obvious inequality:
\begin{equation} \label{trig}
a\cos\phi+b\sin\phi\geq-\sqrt{a^2+b^2}.
\end{equation}
Consider 
$$
g_n(\alpha):=\frac{f_n(\alpha)}{\sin\alpha} = (3n^2-4)+(n^2+4)\cos n\alpha -4n\cot\alpha\sin n\alpha.
$$
Using (\ref{trig}) with $\phi=n\alpha$, we have:
$$
g_n(\alpha) \geq (3n^2-4) - \sqrt{(n^2+4)^2+16n^2 \cot^2\alpha}.
$$
Thus, $g_n(\alpha)>0$ if $(3n^2-4)^2>(n^2+4)^2+16n^2 \cot^2\alpha$, or 
\begin{equation} \label{cot}
|\cot\alpha|<\sqrt{(n^2-4)/2}.
\end{equation}

Secondly, we show that $f_n(\alpha)$ increases when $\alpha$ is sufficiently close to $0$. Namely, let
$$
h_n(\alpha):=\frac{f_n'(\alpha)}{\cos\alpha} = (3n^2-4)(1-\cos n\alpha)-n^3\tan\alpha\sin n\alpha.
$$
Using some trigonometry, we see that $h_n(\alpha)>0$ if $(3n^2-4) \tan (n\alpha/2)>n^3 \tan\alpha$, or
\begin{equation} \label{trigineq}
\frac{\tan (n\alpha/2)}{n\tan (\alpha/2)}>\frac{n^2}{3n^2-4}\ \frac{\tan\alpha}{\tan (\alpha/2)}.
\end{equation}
We want to show that (\ref{trigineq}) holds for $\alpha<\pi/(2n)$. To this end, we use the inequality
$$
\tan n\alpha > n\tan\alpha,
$$
that can be easily proved by induction on $n$. Hence, the left hand side of (\ref{trigineq}) is greater than 1. On the other hand, the right hand side of (\ref{trigineq}) is less than 1 for all $n\geq 4$ and $\alpha<\pi/(2n)$. Thus (\ref{trigineq}) holds.

Finally, we need to show that the two above considered cases cover the whole interval of values of $\alpha$, that is, in view of (\ref{cot}), that 
\begin{equation} \label{last}
\cot \frac{\pi}{2n} < \sqrt{\frac{n^2-4}{2}}
\end{equation}
for all $n\geq 4$. Indeed, the ratio of the right and left side of (\ref{last}) increases with $n$, and for $n=4$, this ratio equals $1.01461...$
\proofend

\begin{remark} \label{crcond} 
Critical curves of $I(\alpha)$.
{\rm It is interesting to describe critical curves of the functional $I(\alpha)$. Let $\gamma(t)$ be a $2T$-periodic centrally symmetric curve, parameterized so that the unit Wronskian condition (\ref{Wron}) holds. Let $\gamma_{\pm}(t)=\gamma(t\pm\alpha)$. Then $\gamma$ is critical for $I(\alpha)$ in the class of curves satisfying  (\ref{Wron}) if and only if 
\begin{equation} \label{cond}
3[\gamma'(t),\gamma_+(t)-\gamma_-(t)]+[\gamma(t),\gamma_+'(t)-\gamma_-'(t)]=0
\end{equation}
for all $t$.
We do not dwell on the proof, but let us mention that the infinitesimal perturbations of a curve, preserving the unit Wronskian condition, are given by vector fields of the form
$$
v(t)=f'(t)\gamma(t)-2f(t)\gamma'(t)
$$
where $f(t)$ is an arbitrary smooth function satisfying $f(t+T)=-f(t)$.

Equation (\ref{cond}) holds for central ellipses (each of the two cross-products vanishes), but we do not know whether central ellipses are the only curves satisfying this equation.
}
\end{remark}

\section{Appendix: Blaschke-Santalo inequality and outer billiards} \label{outerbill}

\begin{figure}[hbtp]
\centering
\includegraphics[width=2in]{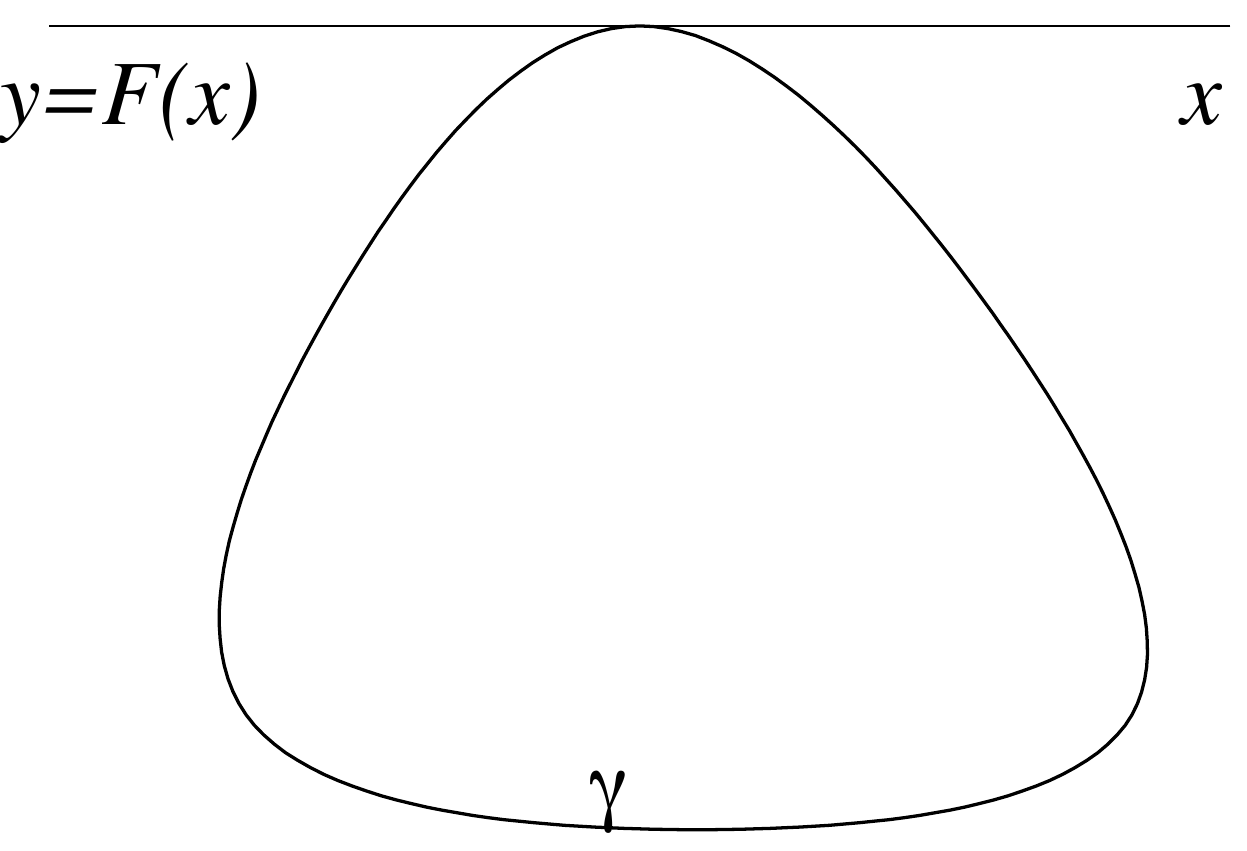}
\caption{Outer billiard map}
\label{outer}
\end{figure}

{\it Outer billiards} (a.k.a. dual billiards) is a discrete time dynamical system in the exterior of a planar convex domain (outer billiard table) defined by the following geometric construction. Let $\gamma$ be the oriented outer billiard curve, the boundary of the outer billiard table, and let $x$ be a point in its exterior. Draw the tangent ray to $\gamma$ from $x$, whose orientation agrees with that of $\gamma$, and reflect $x$ in the tangency point to obtain a new point $y$. The map $F:x\mapsto y$ is the outer billiard transformation, see Figure \ref{outer}. The map $F$ can be defined for convex polygons as well (its domain is then an open dense subset of the exterior of the polygon). See the article \cite{D-T} or the respective chapters of the books \cite{Ta0,Ta3} for a survey of outer billiards. The monograph \cite{Sch3} provides a profound study of outer billiards on a class of quadrilaterals called kites. 

It was observed a long time ago that, after rescaling, the dynamics of the second iteration of the outer billiard map very far away from the outer billiard table is approximated by a continuous motion whose trajectories are closed centrally symmetric curves and which satisfies the second Kepler law:  the area swept by the position vector of a point depends  linearly on time. Without going into details that can be found in \cite{Ta1}, here is an explanation of this phenomenon. 

\begin{figure}[hbtp]
\centering
\includegraphics[width=4.5in]{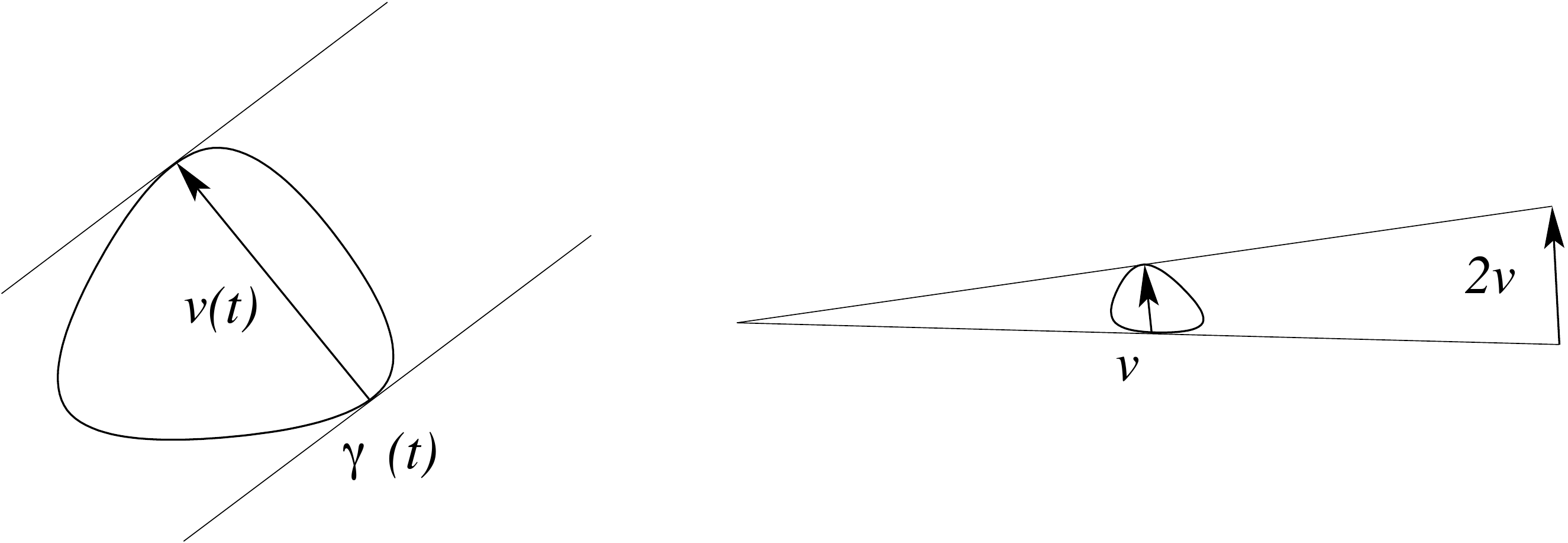}
\caption{Outer billiard map far away from the table}
\label{faraway}
\end{figure}

Let $\gamma(t)$ be the outer billiard curve which we assume to be smooth and strictly convex. Consider the tangent line to $\gamma(t)$. There is another tangent line, parallel to that at $\gamma(t)$;  let $v(t)$ be the vector that connects the tangency  points of the former and the latter. For points $x$ at great distance from $\gamma$ and seen in the direction of $\gamma'(t)$ from $\gamma$, the vector $\overrightarrow{x\ F^2(x)}$ is almost equal to $2v(t)$, see Figure \ref{faraway}. We construct a homogeneous field of directions in the plane: along the ray generated by the vector $\gamma'(t)$, the direction of the field is that of the vector $v(t)$. The trajectories of the second iteration of the outer billiard map ``at infinity" follow the integral curves of this field of directions. These integral curves are all similar; we denote them by $\Gamma$. 

A similar analysis can be performed when the outer billiard curve if a convex polygon, see, e.g., \cite{Sch3}. For example,  if $\gamma$ is a triangle then $\Gamma$ is an affine-regular hexagon. Another example: if $\gamma$ is a curve of constant width then $\Gamma$ is a circle.  If $\gamma$ is a semi-circle then $\Gamma$ is a curve made of two symmetric arcs of orthogonal parabolas. 

In fact, one can describe the curves $\Gamma$ explicitly. Let us assume first that $\gamma$ is centrally symmetric. Then $v(t)=-2\gamma(t)$. 

\begin{lemma} \label{dualcurve}
The curve
\begin{equation} \label{infcurve}
\Gamma(t)=\frac{\gamma'(t)}{[\gamma(t),\gamma'(t)]}
\end{equation}
 is the integral curve of the above defined field of directions. 
\end{lemma}

\proof
Clearly,  $\Gamma(t)$ has the direction of $\gamma'(t)$, and we need to check that $\Gamma'(t)$ is collinear with $\gamma(t)$. Indeed,
$$
\Gamma'=\frac{\gamma''}{[\gamma,\gamma']} - \frac{\gamma' [\gamma,\gamma'']}{[\gamma,\gamma']^2},
$$
hence $[\gamma,\Gamma']=0$.
\proofend

The outer billiard motion ``at infinity" goes along the curve (\ref{infcurve})  with the velocity vector at point $\Gamma(t)$ being equal to $-2\gamma(t)$. This implies that 
\begin{equation} \label{rate}
[\Gamma(t),-2\gamma(t)]=2, 
\end{equation}
which explains Kepler's law. Furthermore, the curve (\ref{infcurve}) is polar dual to $\gamma$: equations (\ref{duality}) hold with $\Gamma=\gamma^*$ (providing another proof to Lemma \ref{dualcurve}).

If $\gamma$ is not centrally symmetric then $\Gamma$ is polar dual to the central symmetrization of $\gamma$, see \cite{Ta1}. The later curve, which we denote by $\bar\gamma$, is the Minkowski half-sum of $\gamma$ and $-\gamma$, its reflection in the origin. In other words, the support function of $\bar\gamma$ is given by the formula
$$
\bar p(t)=\frac{p(t)+p(t+\pi)}{2},
$$
where $p(t)$ is the support function of $\gamma$. The curve $\bar\gamma$ is centrally symmetric and its width in every direction coincides with that of $\gamma$. Of course, if $\gamma$ is centrally symmetric then $\bar\gamma=\gamma$.

The trajectories at infinity $\Gamma$ are defined only up to dilation. Fix one such curve, and let $T$ be the time it takes to traverse the curve moving with the velocity $v$. Scaling the curve $\Gamma$  by some factor,  results in scaling $T$ by the same factor. One can also scale the outer billiard curve $\gamma$: this results in scaling the speed by the same factor and the time by its reciprocal.  To make the time scaling-independent, one multiplies $T$ by  $\sqrt{A(\bar\gamma)/A(\Gamma)}$
where, as before, $A$ denotes the area bounded by a curve. Let us call the result of this scaling of $T$ the {\it absolute time}, and denote it by ${\cal T}$.  

\begin{theorem} \label{abstime}
For any outer billiard curve $\gamma$, the absolute time satisfies
$$
\sqrt{2} \leq {\cal T}\leq \frac{\pi}{2}.
$$
The upper bound is attained only for curves of constant width and their affine images; the lower bound is attained only 
for parallelograms. If $\gamma$ is a centrally symmetric $2n$-gon  then 
$$
{\cal T}\leq n \sin \frac{\pi}{2n},
$$
with equality only for affine-regular $2n$-gons; the same inequality holds for arbitrary $n$-gons.
\end{theorem}

\proof
Let $\Gamma(t)$ be as in Lemma \ref{dualcurve}. According to (\ref{rate}), the rate of change of sectorial area is 2, so the time $T$ equals $(1/2)A(\Gamma)$. Hence ${\cal T}=(1/2)\sqrt{A(\bar\gamma)A(\Gamma)}$. 

By the Blaschke-Santalo inequality, ${\cal T}\leq \pi/2$, with equality only if $\bar\gamma$ is a central ellipse, that is, if $\bar\gamma$ is affine equivalent to a circle. But $\bar\gamma$ is a circle if and only if $\gamma$ has constant width. 

By Mahler's theorem, see \cite{Mah,Lut1},  $\sqrt{2}\leq {\cal T}$, with equality only if $\bar\gamma$ is a parallelogram. This happens if and only if  $\gamma$ is a parallelogram as well. 

If $\gamma$ is a centrally symmetric $2n$-gon, the upper bound follows from Theorem \ref{BS}. Finally, if $\gamma$ is an $n$-gon then $\bar\gamma$ is a centrally symmetric $2n$-gon (it is possible that $\bar\gamma$ has fewer than $2n$ sides but this does not affect the inequality). 
\proofend

\begin{remark}
{\rm It is interesting to mention that outer billiards also ``solve" the isoperimetric problem in Minkoswki geometry. Let a centrally symmetric outer billiard curve  $\gamma$ be the unit circle of  planar Minkowki geometry. Then the trajectory at infinity $\Gamma$ is the unique solution to the isoperimetric problem in this Minkowski geometry: according to Busemann's theorem  \cite{Bu}, the Minkowski length of (a homothetic copy of) $\Gamma$ is minimal among the curves bounding a fixed area.
}
\end{remark}

\bigskip

{\bf Acknowledgments}. It is a pleasure to thank  J. C. Alvarez, Yu. Burago, M. Ghomi, M. Levi, E. Lutwak, V. Ovsienko, I. Pak and especially R. Schwartz for comments and suggestions. The paper was written during my visit at Brown University; I am  grateful to the Department of Mathematics for its hospitality.

\end{document}